\documentclass{SCIYA2017enOL}
\online
\numberwithin{equation}{section}

\newcommand{\R}{\mathbb{R}}

\newcommand{\eps}{ \varepsilon}
\newcommand{\Det}{\mathcal{D}\textrm{et }}

\usepackage{epsfig}
\usepackage{hyperref}
\usepackage{epstopdf}

\begin{document}

\ensubject{fdsfd}

\ArticleType{ARTICLES}
\SubTitle{Dedicated to Professor Jean-Yves Chemin on the Occasion of His {\rm 60}th Birthday}
\Year{}
\Month{}%
\Vol{}
\No{}
\BeginPage{1} %
\DOI{}
\ReceiveDate{December 6, 2018}
\AcceptDate{March 14, 2019}

\title[Monge-Amp\'ere equation]{Very weak solutions to the two dimensional Monge-Amp\'ere equation}

\author[1]{Wentao Cao}{wentao.cao@math.uni-leipzig.de,}
\author[1$\ast$]{L\'aszl\'o Sz\'ekelyhidi Jr.}{laszlo.szekelyhidi@math.uni-leipzig.de}

\AuthorMark{W. Cao}

\AuthorCitation{W. Cao, L. Sz\'ekelyhidi}

\address[1]{Institut f\"{u}r mathematik, Universit\"{a}t Leipzig, D-04109, Leipzig, Germany}

\abstract{In this short note we revisit the convex integration approach to constructing very weak solutions to the 2D Monge-Amp\'ere equation with H\"older-continuous first derivatives of exponent $\beta<1/5$. Our approach is based on combining the approach of Lewicka-Pakzad \cite{Lewicka:2015gz} with a new diagonalization procedure which avoids the use of conformal coordinates, which was introduced by the second author with De Lellis and Inauen in \cite{DeLellis:2015wm} for the isometric immersion problem.}



\maketitle

\section{Introduction}

 In this short note we consider very weak solutions of the 2D Monge-Amp\'ere equation. As pointed out by T.~Iwaniec in \cite{Iwaniec:2001uw}, the Hessian in two real variables can be written in various weak forms. In particular, using the identity
$$
\partial_{11}v\partial_{22}v-(\partial_{12}v)^2=
\partial_{12}(\partial_1v\partial_2v)-\tfrac{1}{2}\partial_{22}(\partial_1v)^2-\tfrac{1}{2}\partial_{11}(\partial_2v)^2
$$
one can define the \emph{very weak Hessian} for $v\in W^{1,2}_{loc}$. Denoting the right hand side of the above formula (in vector notation) as
$$
\Det\nabla^2v:=-\frac{1}{2}\textrm{curl }\textrm{curl}(\nabla v\otimes\nabla v).
$$
we consider \emph{very weak solutions} of the Monge-Amp\'ere equation in the sense
\begin{equation}\label{e:mongeampere}
\Det\nabla^2v=f \quad \text{ in } \mathcal{D}'(\Omega),
\end{equation}
i.e. in the sense of distributions, where $\Omega\subset\R^2$ is a simply connected open subset.

It was noted by M.~Lewicka and M.R.~Pakzad in \cite{Lewicka:2015gz} that there is a close connection between $C^1$ solutions of \eqref{e:mongeampere} and $C^1$ isometric immersions $\Omega\subset\R^2\to\R^3$. In particular, an adaptation of convex integration, following the famous Nash-Kuiper technique \cite{Nash:1954vt,Kuiper:1955th} and its extension to $C^{1,\beta}$ solutions in \cite{Conti:2012wl}, can be used to construct $C^{1,\beta}$ solutions to \eqref{e:mongeampere} (along with analogous h-principle statements). In particular, the threshold regularity that can be reached is $\beta<1/7$, as in \cite{Borisov:1965wf,Conti:2012wl}. The authors in \cite{Lewicka:2015gz} also prove a rigidity statement for $\beta>2/3$ for the degenerate case $f=0$, following work in \cite{Jerrard:2009go,Kirchheim:2003wp,Pakzad:2004vf}. Furthermore, in \cite{Lewicka:2015gz} it is noted as a footnote on p2 that the recent idea in \cite{DeLellis:2015wm} of using a transformation to conformal coordinates could potentially be used in the present setting as well to increase the threshold regularity to $\beta<1/5$. Our purpose in this short note is to give a much simpler proof of this extension, which avoids use of conformal coordinates.

Why are we interested in this problem?
A recurring theme in the modern theory of nonlinear partial differential equation (PDE) is that, whilst smooth solutions of a certain equation automatically satisfy a \emph{derived} equation leading to additional properties such as uniqueness and higher regularity, weak or distributional solutions may exist for which the derived equation does not hold. Invariably, in such situations an important question is the precise regularity threshold, which guarantees validity of the derived PDE. Probably the most well-known example is the linear transport equation: if $\rho$ is a solution to $\partial_t\rho+v\cdot\nabla\rho=0$ with sufficiently regular velocity $v$, any composition of the form $\beta\circ \rho$ with arbitrary smooth functions $\beta$ is also a solution. This observation and the associated regularity threshold underlies the powerful theory of renormalized solutions \cite{DiPerna:1989vo,Depauw:2003wl,Ambrosio:2017hb,Modena:2017tp}. Further well-known examples arise in various PDE in fluid mechanics, starting with the inviscid Burger's equation, entropy solutions of hyperbolic conservation laws and the incompressible Euler equations, in connection with turbulence and anomalous dissipation \cite{Constantin:1994vn,Isett:2016to,Buckmaster:2017uza}. Finally, we mention the system of equations describing isometric immersions of surfaces in Euclidean space, which is closest to our focus in this paper: whereas for classical isometric immersions Gauss' Theorema Egregium holds and plays an important role in the proof of rigidity of smooth convex surfaces \cite{Herglotz:1943je}, the Nash-Kuiper theorem \cite{Nash:1954vt,Kuiper:1955th} and its recent extensions \cite{Borisov:1965wf,Conti:2012wl,DeLellis:2015wm} produce (weak, merely $C^1$) solutions of the isometric embedding problem to which no reasonable notion of curvature can be associated.

In certain cases the original PDE under consideration is already the ``derived'' equation, for which very weak solutions may be defined via a kind of primitive: indeed, often in fluid mechanics the equations of motion are first formulated in vorticity form (e.g.~for 2D incompressible Euler) and the velocity formulation is treated as the primitive form - hence also the terminology used for the primitive equations of geophysical flows \cite{Lions:1992fs}. Another recent example belonging to this class is the inviscid SQG equation, where weak solutions enjoy a weak form of rigidity (more precisely, stability under weak convergence) but an entirely different class of ``very weak'' solutions may be constructed via convex integration \cite{Resnick:1995tr,Buckmaster:2016vy}. The present context belongs to this class of problems - thus, in light of rigidity results for the Monge-Amp\'ere equation in the Sobolev setting \cite{Pakzad:2004vf,Sverak:2000wa}   an interesting future direction is to understand the threshold between rigidity versus flexibility in a Sobolev or Besov scale of spaces.

Returning to problem \eqref{e:mongeampere}, we recall that
the kernel of the differential operator \textit{curl }\textit{curl} in 2D consists of fields of the form $(\textrm{sym}\nabla w)\overset{\textrm{def}}{=}\frac12(\nabla w+\nabla w^T)$ on a simply connected domain $\Omega$. Hence $v$ is a solution of the homogeneous case $f=0$ in \eqref{e:mongeampere} if and only there exists a mapping $w:\Omega\rightarrow\R^2$ such that
\begin{equation}\label{e:homo}
\frac{1}{2}\nabla v\otimes\nabla v+\textrm{sym }\nabla w=0  \text{ in } \Omega.
\end{equation}
The non-homogeneous case is then equivalent to the problem of finding a pair $(v, w)$ such that
\begin{equation}
\label{e:matrix-reduced}
\frac{1}{2}(\nabla v\otimes\nabla v)+\textrm{sym }\nabla w=A
\end{equation}
together with the auxiliary linear problem
\begin{equation}\label{e:curlcurl}
-\textrm{curl }\textrm{curl } A=f.
\end{equation}
In \cite{Lewicka:2015gz} the authors solve \eqref{e:matrix-reduced} by adapting the convex integration technique of Nash from \cite{Nash:1954vt} as developed for $C^{1,\beta}$ solutions in \cite{Conti:2012wl}. In this note we show that the corrector field $\textrm{sym }\nabla w$ can be used to diagonalize the deficit matrix directly, without resorting to conformal coordinates, and in this way we can directly apply the Nash convex integration scheme with 2 steps rather than 3 -   this is responsible for the improvement $1/7$ to $1/5$ (see \cite{SzekelyhidiJr:2014tu} for an exposition explaining this point). Our main result is

\begin{theorem}\label{t:c1.2weak}
For any $f\in L^p$ with $p>\frac{5}{4}$ on an open bounded simply connected domain $\Omega\subset\R^2$ with $C^{1, 1}$ boundary and any $\beta<\frac{1}{5},$  $C^{1, \beta}$ weak solutions to \eqref{e:mongeampere} are dense in $C^0(\bar\Omega).$ More precisely, for any ${v}^\flat\in C^{0}(\bar\Omega),$ for any $\eps>0,$ there exists a  weak solution $v\in C^{1, \beta}$ to the Monge-Amp\'ere equation, i.e.
$$\Det \nabla^2v=f \quad \text{ in } \Omega,$$
such that
$$\|v-v^\flat\|_0\leq\eps.$$
\end{theorem}

The connection between very weak solutions of \eqref{e:mongeampere} and $C^{1,\beta}$ isometric immersions of course runs much deeper than merely availability of the same method. Just as it is well-known that via the graphical representation of surfaces, the (classical, $C^2$) isometric immersion problem is equivalent to the Darboux equation, which is of Monge-Amp\'ere type (see for instance \cite{Han:2006fo}), the present very weak formulation \eqref{e:mongeampere} corresponds to a first order approximation of the (weak) isometric embedding problem. Indeed, following \cite{Lewicka:2015gz}, consider a $1$-parameter family of deformations by an out-of-plane displacement $v:\Omega\to\R$ and an in-plane (quadraticly scaled) displacement $w:\Omega\to\R^3$,
\begin{equation}\label{e:approximation}
\phi_t=id+ tve_3+t^2w: \Omega\rightarrow\R^3,
\end{equation}
with $e_3=(0, 0, 1)$. Then  \eqref{e:homo} can be seen as an  equivalent condition for deformations $\phi_t$ to form a second order infinitesimal bending, i.e. to make
\begin{align*}
\nabla\phi_t^T\nabla\phi_t=\textrm{Id}+o(t^2),
\end{align*}
where $\textrm{Id}$ stands for 2-dimensional identity matrix.

To close this introduction we wish to emphasize that our result (as well as the result in \cite{Lewicka:2015gz}) does not solve the Dirichlet problem, where one would couple \eqref{e:mongeampere} with a boundary condition $v=\phi$ on $\partial\Omega$. This would be an interesting problem in its own right: indeed, we know that the full Dirichlet condition for the isometric immersion problem is over-determined, even in the weak setting \cite{Hungerbuhler:2017hn,Cao:2018uf}, but prescribing $v$ on the boundary corresponds to partial Dirichlet data (cf. \eqref{e:approximation}), and we believe that an analogous result to Theorem \ref{t:c1.2weak} should be true.

\section{Preliminaries}\label{s:preliminary}
We introduce some notation, function spaces and basic lemmas in this section. For any multi-index $\beta$  and for any function $h:\Omega\rightarrow\R$, we define the supreme norm as
\begin{equation*}
\|h\|_0=\sup_{\Omega}h, ~~\|h\|_m=\sum_{j=0}^m\max_{|\beta|=j}\|\partial^{\beta} h\|_0,
\end{equation*}
and H\"older semi norms as
\begin{align*}
[h]_{\alpha}&=\sup_{x\neq y}\frac{|h(x)-h(y)|}{|x-y|^\alpha},\\
[h]_{m+\alpha}&=\max_{|\beta|=m}\sup_{x\neq y}\frac{ |\partial^{\beta} h(x)-\partial^{\beta} h(y)|}{|x-y|^\alpha},
\end{align*}
for any $0<\alpha\leq 1$.
The H\"older norms are given as
$$
\|h\|_{m+\alpha}=\|h\|_m+[h]_{m+\alpha}.
$$
We recall the interpolation inequality for H\"older norms
$$
[h]_r\leq C\|h\|_0^{1-\frac{r}{s}}[h]_s^{\frac{r}{s}}
$$
for $s>r\geq0.$  In particular,
\begin{equation}\label{e:holder-alpha}
\|h\|_\alpha\leq \|h\|_0+2\|h\|_0^{1-\alpha}[h]_1^\alpha \text{ for any } \alpha\in[0, 1].
\end{equation}

We also recall some estimates of the regularization of H\"older functions.
\begin{lemma}\label{l:mollification}
For any $r, s\geq0,$ and $0<\alpha\leq1,$ we have
\begin{align*}
&[f*\varphi_l]_{r+s}\leq Cl^{-s}[f]_r,\\
&\|f-f*\varphi_l\|_r\leq Cl^{1-r}[f]_1, \text{ if } 0\leq r\leq1,\\
&\|(fg)*\varphi_l-(f*\varphi_l)(g*\varphi_l)\|_r\leq Cl^{2\alpha-r}\|f\|_\alpha\|g\|_\alpha,
\end{align*}
with constant $C$ depending only on $s, r, \alpha, \varphi.$
\end{lemma}
Other properties about H\"older norm can be found in references such as \cite{Conti:2012wl,DeLellis:2015wm,Buckmaster:2017uza}. The norms of a $2\times 2$  matrix $P$ is defined as
\begin{align*}
|P|:=\sup_{\xi\in S^{n-1}}|P\xi|.
\end{align*}
In the paper, $C(\cdot)$ denotes constants depending the parameter in the bracket. We also recall the following lemma about the corrugation functions for the two dimensional Monge-Amp\'ere equation in \cite{Lewicka:2015gz}.
\begin{lemma}\label{l:Gamma}
There exists a smooth 1-periodic field $\Gamma=(\Gamma_1, \Gamma_2)\in C^{\infty}([0, \infty)\times\R, \R^2)$ such that the following holds for any $(s, t)\in[0, \infty)\times\R:$
\begin{equation}\label{e:gamma12}
\Gamma(s, t+1)=\Gamma(s, t), \qquad \frac{1}{2}|\partial_t\Gamma_1(s, t)|^2+\partial_t\Gamma_2(s, t)=s^2,
\end{equation}
along with the following estimates: for any nonnegative integer $k,$
\begin{align*}
&|\partial_t^k\Gamma_1(s, t)|+|\partial_s\partial_t^k\Gamma_2(s, t)|\leq Cs, \\
&|\partial_s\partial_t^k\Gamma_1(s, t)|\leq C, \quad |\partial_t^k\Gamma_2(s, t)|\leq Cs^2.
\end{align*}
\end{lemma}
Indeed, one choice of such $\Gamma$ in \cite{Lewicka:2015gz} is
\begin{equation}\label{e:gamma}
\Gamma_1(s, t)=\frac{s}{\pi}\sin (2\pi t),\qquad  \Gamma_2(s, t)=-\frac{s^2}{4\pi}\sin(4\pi t).
\end{equation}

\section{Diagonalisation}\label{s:diagonal}
In this section, we show that diagonalization of the deficit matrix can be achieved by solving a planar div-curl system for the corrector field is equivalent $\textrm{sym}\nabla w$.
\begin{proposition}\label{p:diag}
For any $j\in\mathbb{N}, 0<\alpha<1,$ there exist constants $M_1, M_2, \sigma_1$ depending only on $j, \alpha$ such that the following statement hold. If $D\in C^{j, \alpha}(\Omega, \R_{sym}^{2\times2})$ satisfies
\begin{equation}
\label{e:sigma}
\|D-\textrm{Id}\|_\alpha\leq \sigma_1,
\end{equation}
then there exists $\Phi\in C^{j+1, \alpha}(\Omega, \R^2)$ and $d\in C^{j, \alpha}(\Omega, \R)$ such that
\begin{equation}
\label{e:diag}
D+\textrm{sym }\nabla\Phi=d^2\textrm{Id}
\end{equation}
and the following estimates hold:
\begin{align}
&\|d-1\|_\alpha+\|\nabla\Phi\|_\alpha\leq M_1\|D-\textrm{Id}\|_{\alpha}; \label{e:d-estimate-alpha}\\
&[d]_{j,\alpha}+[\nabla\Phi]_{j,\alpha}\leq M_2\|D-\textrm{Id}\|_{j+\alpha}. \label{e:d-estimate-high}
\end{align}
\end{proposition}
\begin{proof}
Writing $D=\tfrac{1}{2}(\textrm{tr}D)\textrm{Id}+\mathring{D}$, where $\mathring{D}$ is the traceless part of $D$, from the assumption \eqref{e:sigma} we deduce
\begin{equation} \label{e:dijestimate-1}
\|\mathring{D}\|_\alpha\leq \sigma_1, \quad  \tfrac12\textrm{tr}D\geq 1-\sigma_1.
\end{equation}
Note that \eqref{e:diag} amounts to the system
\begin{equation}\label{e:dijequation-1}
\begin{split}
\displaystyle &D_{12}+\frac{1}{2}(\partial_{x_2}\Phi_1+\partial_{x_1}\Phi_2)=0,\\
\displaystyle &D_{11}+\partial_{x_1}\Phi_1=D_{22}+\partial_{x_2}\Phi_2(x)=d^2.
\end{split}
\end{equation}
Equivalently \eqref{e:dijequation-1} can be written as
\begin{equation} \label{e:dijequation-2}
\begin{split}
\displaystyle&\partial_{x_2}\Phi_1+\partial_{x_1}\Phi_2=\textrm{curl}_{x_1, x_2}(-\Phi_1, \Phi_2)=-2D_{12},\\
\displaystyle&\partial_{x_2}\Phi_2-\partial_{x_1}\Phi_1=\textrm{div}_{x_1, x_2}(-\Phi_1, \Phi_2)=D_{11}-D_{22},
\end{split}
\end{equation}
which is a planar div-curl system for vector function $(-\Phi_1, \Phi_2)$. Setting $\Phi_1=-\partial_{x_1}\varphi-\partial_{x_2}\psi$ and $\Phi_2=\partial_{x_2}\varphi-\partial_{x_1}\psi$ leads to the equations
$$
\Delta \varphi=D_{11}-D_{22},\quad \Delta\psi=2D_{12}.
$$
Thus, solving the associated Dirichlet problems (i.e.~with $\varphi=\psi=0$ on $\partial\Omega$) and using standard Schauder estimates, we obtain a solution $\Phi\in C^{j+1, \alpha}(\Omega, \R^2)$ of \eqref{e:dijequation-2} with
$$
\|\nabla\Phi\|_{j+\alpha}\leq C_{j,\alpha}\|\mathring{D}\|_{j,\alpha}.
$$
Next, from \eqref{e:dijequation-1} we see that
$$
d^2=\tfrac12\textrm{tr}D+\tfrac12\textrm{div}\Phi\geq 1-(1+C_{\alpha})\sigma_1,
$$
so that, choosing $\sigma_1$ sufficiently small ensures that $d^2\geq 1/2$, from which we easily deduce the required estimates in \eqref{e:d-estimate-alpha}-\eqref{e:d-estimate-high}.
\end{proof}

\section{Proof of Theorem \ref{t:c1.2weak}}\label{s:proof}
In this section, we will prove Theorem \ref{t:c1.2weak} by iteratively  constructing a sequence of  subsolutions $(v_q, w_q).$ Our proof is also divided into four steps.

\textbf{Step 1. Formulation.}
As in \cite{Lewicka:2015gz}, the problem of seeking $C^{1, \beta}$ solution can be formulated to construct solutions of the following equations
\begin{align}
&-\textrm{curl }\textrm{curl }A=f,\label{e:A}\\
 &A=\frac{1}{2}\nabla v\otimes\nabla v+\textrm{sym}\nabla w.\label{e:matrix}
\end{align}
For any $f\in L^p(\Omega),$ the Dirichlet problem
\begin{align*}
-\triangle u=f \text{ in }\Omega,  \quad u=0 \text{ on }\partial\Omega
\end{align*}
admits a $W^{2, p}$ solution $u$, then Morrey's theorem further implies that $u\in C^{0, \kappa}(\bar\Omega)$ for $\kappa=2-\frac{2}{p}.$ Set $A=(u+c)\textrm{Id}$ where $c$ is a constant to be fixed, then
\begin{align*}
-\textrm{curl }\textrm{curl }A=-\triangle(u+c)=f.
\end{align*}
and $A\in C^{0, \kappa}(\Omega).$ So we get a weak solution to \eqref{e:A}. Since $C^\infty(\bar\Omega)$ is dense in $C^0(\bar\Omega),$ we can assume $v^\flat\in C^{\infty}(\bar\Omega)$. Take $ w^\flat=0,$ then choose $c$ large enough to make
$$D^\flat=A-\frac{1}{2}\nabla v^\flat\otimes\nabla v^\flat-\textrm{sym }\nabla w^\flat\geq 2\bar\delta\textrm{Id}$$
with some $\bar\delta>0$. It is only remained to solve matrix equation
\eqref{e:matrix}. We will construct $v, w$ satisfying \eqref{e:matrix} through convex integration.

\textbf{Step 2. Initial approximation.} We can assume $A$ is smooth in this step since the mollification error in this step will not be iterated.  Due to
\begin{align*}
A-\frac{1}{2}\nabla v^\flat\otimes\nabla v^\flat-\textrm{sym }\nabla w^\flat-\bar\delta\textrm{Id}\geq\bar\delta\textrm{Id},
\end{align*}
using Lemma 1 in \cite{Nash:1954vt} (see also Lemma 1.9 in \cite{SzekelyhidiJr:2014tu}), we obtain the decomposition
\begin{align*}
A-\frac{1}{2}\nabla v^\flat\otimes\nabla v^\flat-\textrm{sym }\nabla w^\flat-\bar\delta\textrm{Id}=\sum_{i=1}^Na_i^2\nu_i\otimes\nu_i
\end{align*}
for some $\nu_k\in S^{1}, a_k\in C^{\infty}(\bar\Omega)$ and some integer $N.$ Define iteratively the smooth mappings $\bar v_0=v^\flat, \bar v_1,\cdots, \bar v_N,$ and  $\bar w_0=w^\flat, \bar w_1,\cdots, \bar w_N,$ as in \cite{Lewicka:2015gz},
\begin{align*}
\bar v_i&=\bar v_{i-1}+\frac{1}{\mu_i}\Gamma_1(d(x), \mu_i x\cdot\nu_i),\\
\bar w_i&=\bar w_{i-1}-\frac{1}{\mu_i}\Gamma_1(d(x), \mu_i x\cdot\nu_i)\nabla\bar v_{i-1}+\frac{1}{\mu_i}\Gamma_2(d(x), \mu_i x\cdot\nu_i)\nu_i,
\end{align*}
where  frequencies $1\leq\mu_1\leq\mu_2\leq\cdots\leq\mu_N$ are to be fixed. Denote
\begin{equation*}
\mathcal{E}_i=\left(\frac{1}{2}\nabla\bar{v}_i\otimes\nabla\bar{v}_i+\textrm{sym}\nabla\bar{w}_i\right)-\left(a_i^2 \nu_i\otimes \nu_i+\frac{1}{2}\nabla\bar{v}_{i-1}\otimes\nabla\bar{v}_{i-1}+\textrm{sym}\nabla\bar{w}_{i-1}\right),
\end{equation*}
then
\begin{align*}
A-\frac{1}{2}\nabla\bar v_N\otimes\nabla \bar v_N-\textrm{sym }\nabla \bar w_N=\bar\delta\textrm{Id}-\sum_{i=1}^N\mathcal{E}_i.
\end{align*}
Similar to the calculation in Step 2 of the proof of Proposition \ref{p:stage}, it is not hard to get
\begin{align*}
\|\mathcal{E}_1\|_0\leq\frac{C(v^\flat, w^\flat)}{\mu_1}, \quad \|\mathcal{E}_1\|_1\leq C(v^\flat, w^\flat), \quad \|\bar{v}_1-v^\flat\|_0+\|\bar{w}_1-w^\flat\|_0\leq \frac{C}{\mu_1}.
\end{align*}
Hence interpolation of H\"older spaces gives
\begin{align*}
\|\mathcal{E}_1\|_\alpha\leq\frac{\sigma_0}{2N}\bar\delta, \quad \|\bar{v}_1-v^\flat\|_0+\|\bar{w}_1-w^\flat\|_0\leq \frac{\eps}{2N},
\end{align*}
provided choosing $\mu_1$ large.  Analogously, for $\bar{v}_2$, we have
\begin{align*}
\|\mathcal{E}_2\|_0\leq\frac{C(v^\flat, w^\flat, \mu_1)}{\mu_2}, \quad \|\mathcal{E}_2\|_1\leq C(v^\flat, w^\flat, \mu_1), \quad \|\bar{v}_2-\bar{v}_1\|_0+\|\bar{w}_2-\bar{w}_1\|_0\leq \frac{C}{\mu_2}.
\end{align*}
Again using interpolation and taking $\mu_2$ large we gain
\begin{align*}
\|\mathcal{E}_2\|_\alpha\leq\frac{\sigma_0}{2N}\bar\delta, \quad \|\bar{v}_2-\bar{v}_1\|_0+\|\bar{w}_2-\bar{w}_1\|_0\leq \frac{\eps}{2N}.
\end{align*}
Inductively, we can take $\mu_i, i=3, \cdots, N$ such that finally we obtain
\begin{align*}
\|A-\frac{1}{2}\nabla\bar v_N\otimes\nabla \bar v_N-\textrm{sym }\nabla \bar w_N-\bar\delta\textrm{Id}\|_\alpha\leq\sum_{i=1}^N\|\mathcal{E}_i\|_\alpha\leq\frac{\sigma_0}{2}\bar{\delta},
\end{align*}
and
\begin{align*}
\|\bar{v}_N-v^\flat\|_0+\|\bar{w}_N-w^\flat\|_0\leq \frac{\eps}{2}.
\end{align*}

Having constructed a subsolution which satisfies the conditions of the diagonalisation Proposition \ref{p:diag}, we are now ready for the iteration. First we define some parameters which are similar to those in \cite{DeLellis:2015wm}. The  parameter $\alpha$ is a H\"older index which is assumed to be smaller than a geometric constant $\alpha_0$, i.e.
\begin{equation}\label{e:alpha}
0<\alpha<\alpha_0.
\end{equation}
The amplitude parameter  $\delta_q$ and   frequency parameter $\lambda_q$ are defined as
\begin{equation}\label{e:delta-q-lambda-q}
\delta_q=a^{-b^q},\qquad \lambda_q=a^{cb^{q+1}},
\end{equation}
where $q$ is nonnegative integers, $a$ is a large constant  and $b, c>1$ to be prescribed.
The mollification parameter $\ell$ is defined as
\begin{equation} \label{e:l}
\ell^{2-\alpha}=\frac{\delta_{q+1}}{K_1\delta_q\lambda_q^2},
\end{equation}
and the other frequency parameter $\mu$ is
\begin{equation}\label{e:mu}
\mu=\frac{K_2\delta_{q+1}\lambda_{q+1}^\alpha}{\delta_{q+2}\ell},
\end{equation}
where  $K_1, K_2>1$ are suitable large constants to be chosen in the later proof and depending only on $\alpha, A, \sigma_0, K$ but independent of $a.$  Then we have

\begin{lemma}\label{l:b-c-kappa}
There exist $b, c$ such that when $a\gg 1$ is large enough the following holds:
\begin{align}
&\delta_q\lambda_q^2\geq1, \quad \delta_{q+1}<\delta_q\leq1, \quad \lambda_{q+1}>\lambda_q\geq1,\label{e:parameter-monotone}\\
&\lambda_{q+1}^{1-\alpha}\geq \frac{\delta_{q+1}^2\lambda_{q+1}^\alpha}{\delta_{q+2}^2\ell}\geq\mu\geq\ell^{-1}\geq\ell^{-1+\frac{\alpha}{2}}\geq\frac{\delta_q^{1/2}}{\delta_{q+1}^{1/2}}\lambda_q\geq\lambda_q;\label{e:bc-admissible}
\end{align}
and there is $\kappa$ such that
\begin{equation}\label{e:kappa-admissible}
\delta_{q+1}^\kappa\lambda_{q+1}^{\alpha(2-\alpha)}\leq\delta_{q+2}^{2-\alpha}\delta_{q}^\kappa\lambda_q^{2\kappa}.
\end{equation}
\end{lemma}
For simplicity, we also denote
\begin{align*}
A_q=A-\delta_{q+1}\textrm{Id}, \quad D_q=A_q-\frac{1}{2}\nabla v_q\otimes\nabla v_q-\textrm{sym}\nabla w_q.
\end{align*}
\begin{proposition}\label{p:stage}{\bf[Stage]}
There exists a constant $\alpha_0$ such that for any $\alpha$ in \eqref{e:alpha} we can seek positive constants $\sigma_0\leq\frac{\sigma_1}{3}$ with $\sigma_1$ in Proposition \ref{p:diag}, and $K_0$ with the following property. Assume $b, c$ satisfy \eqref{e:bc-condition}, $\kappa$ satisfies \eqref{e:kappa-condition}, $K\geq K_0$ independent of $a$, and $\delta_q, \lambda_q$ defined in \eqref{e:delta-q-lambda-q} with $a$ sufficiently large such that
$$a>a_0(\alpha, b, c, A, K).$$
Let $v_q\in C^2(\bar\Omega, \R)$, $w_q\in C^2(\bar\Omega, \R^2)$ and $A\in C^{0, \kappa}(\bar{\Omega}, \R^{2\times2}_{\textrm{sym}})$ satisfy
\begin{align}
\|D_q\|_\alpha&\leq\sigma_0\delta_{q+1},\label{e:dq}\\
\|v_q\|_{2}+\|w_q\|_2&\leq K\delta_{q}^{1/2}\lambda_q,\label{e:vwq}
\end{align}
then there exist $v_{q+1}\in C^2(\bar\Omega, \R), w_{q+1}\in C^2(\bar\Omega, \R^2)$ such that
\begin{align}
&\|D_{q+1}\|_0\leq\frac{\sigma_0}{3}\delta_{q+2}\lambda_{q+1}^{-\alpha},\quad
\|\nabla D_{q+1}\|_0\leq\frac{\sigma_0}{3}\delta_{q+2}\lambda_{q+1}^{1-\alpha},\label{e:dq+1} \\
&\|v_{q+1}-v_q\|_0\leq \delta_{q+1}^{1/2}\lambda_{q+1}^{-\gamma}, \quad \|w_{q+1}-w_q\|_0\leq M\delta_{q+1},\label{e:vwq+1--0}\\
&\|\nabla(v_{q+1}-v_q)\|_0+\|\nabla(w_{q+1}-w_q)\|_0\leq M\delta_{q+1}^{1/2}, \label{e:vwq+1--1} \\
&\|\nabla^2v_{q+1}\|_0+\|\nabla^2w_{q+1}\|_0\leq K\delta_{q+1}^{1/2}\lambda_{q+1}, \label{e:vwq+1--2}
\end{align}
with some constant $M$ depending only on $\sigma_0$ and $\gamma$ depending only on $\alpha, b, c$.
\end{proposition}
We postpone the proof of Proposition \ref{p:stage}, which corresponds to ''one stage'' in convex integration schemes, into Section \ref{s:iterate-prop}.

 In Step 2 we obtained a first approximation $(\bar{v}_N, \bar{w}_N)$ for  $(v^\flat, w^\flat),$ but as we can see from the assumptions in Proposition \ref{p:stage}, we require the information of the size of the second derivatives of $(\bar{v}_N, \bar{w}_N)$ such that the constant $K$ does not dependent on $a$. Hence we shall use one time Proposition \ref{p:stage} first to get a further approximation $v_0, w_0$ for our later iteration.

\textbf{Step 3. Further approximation.} Mollify $A$ with length-scale $\ell_0$ to get $\tilde{A},$ and set
$$\bar{A}=\tilde A-\bar\delta\textrm{Id}, \quad \bar{D}=\bar{A}-\frac{1}{2}\nabla\bar v_N\otimes\nabla \bar v_N-\textrm{sym }\nabla \bar w_N,$$
then
\begin{align*}
\|\bar{D}\|_\alpha\leq\frac{\sigma_0}{2}\bar{\delta}+C\|A\|_\kappa\ell_0^{\kappa-\alpha}\leq\frac{2\sigma_0}{3}\bar{\delta}, \quad \|\bar{v}_N\|_2+\|\bar{w}_N\|_2\leq C(\mu_N, v^\flat, w^\flat),
\end{align*}
by choosing $\ell_0$ small and $\alpha<\kappa$. First from Step 2 in Proposition \ref{p:stage}, we obtain a  smooth function
$\bar{w}_*$ such that
\begin{align*}
A-\delta_1\textrm{Id}-\frac{1}{2}\nabla\bar v_N\otimes\nabla \bar v_N-\textrm{sym }\nabla \bar w_N=d_0(x)^2(e_1\otimes e_1+e_2\otimes e_2).
\end{align*}
Then define
\begin{align*}
\hat{v}&=\bar v_{N}+\frac{1}{\theta}\Gamma_1(d(x), \theta x\cdot e_1),\\
\hat{w}&=\bar w_{N}+\bar{w}_*-\frac{1}{\theta}\Gamma_1(d(x), \theta x\cdot e_1)\nabla\bar v_{N}+\frac{1}{\theta}\Gamma_2(d(x), \theta x\cdot e_1) e_1,
\end{align*}
and
\begin{align*}
v_0&=\hat{v}+\frac{1}{\lambda}\Gamma_1(d(x), \lambda x\cdot e_2),\\
w_0&=\hat{w}-\frac{1}{\lambda}\Gamma_1(d(x), \lambda x\cdot e_2)\nabla\hat{v}+\frac{1}{\lambda}\Gamma_2(d(x), \lambda x\cdot e_2) e_2.
\end{align*}
Similar computations to Step  3 and Step 4 in Proposition \ref{p:stage} contribute to
\begin{align*}
&\|v_0\|_2+\|w_0\|_2\leq C(\sigma_0, \mu_i, v^\flat, w^\flat)\bar\delta^{1/2}\lambda,\\
&\|A-\frac{1}{2}\nabla v_0\otimes\nabla v_0-\textrm{sym }\nabla w_0-\delta_1\textrm{Id}\|_\alpha\\
&\quad \leq C\|A\|_\kappa\ell_0^{\kappa-\alpha}+ C(\sigma_0, \mu_i, v^\flat, w^\flat)(\bar{\delta}^{1/2}\theta^{2\alpha-1}+\bar\delta\theta\lambda^{\alpha-1}).
\end{align*}
To make sure that $v_0, w_0$ satisfy the assumptions of Proposition \ref{p:stage}, we take
\begin{align*}
\theta=C_1\delta_1^{-1/(1-2\alpha)}, \qquad \lambda=C_2\theta^{1/(1-\alpha)}\delta_1^{-1/(1-\alpha)}
\end{align*}
After choosing the constants $C_1, C_2$ large enough, we can verify \eqref{e:dq} . Then
\begin{align*}
&\|v_0\|_2+\|w_0\|_2\leq C_3\delta_1^{-2/(1-2\alpha)}
\end{align*}
with $C_3$ depending only on $v^\flat, w^\flat, A, \alpha$.  However, we shall show
\begin{align*}
\delta_1^{-2/(1-2\alpha)}\leq \delta_0^{1/2}\lambda_0,
\end{align*}
which after taking logarithms in base $a$ implies
\begin{equation}\label{e:bccondition2}
c\geq\frac{1}{2b}+\frac{2}{1-2\alpha}.
\end{equation}
Hence we also require $b, c$ satisfy \eqref{e:bccondition2}.

\textbf{Step 4. Iteration and conclusion.} Now we are ready to iterate based on Proposition \ref{p:stage}. Fix $\alpha, b, c$ satisfy \eqref{e:bc-condition} and \eqref{e:bccondition2}, then for any large enough $a$ we can construct $(v_0, w_0)$ as in Step 3 such that
\begin{equation}
\label{e:first-difference}
\|v_0-v^\flat\|_0\leq\frac{\eps}{2}
\end{equation}
and all the assumption of Proposition \ref{p:stage} are satisfied for $q=0$ with $K\geq C_3$ independent of $a.$ Thus we can construct  a sequence of approximations $(v_q, w_q)$ by applying Proposition \ref{p:stage}. The sequence satisfies all the conclusions of Proposition \ref{p:stage}. From \eqref{e:vwq+1--1}, it is easy to find that $\{(v_q, w_q)\}$ is a Cauchy sequence in $C^{1}(\bar{\Omega}).$  Let $v, w\in C^{1}(\bar\Omega)$  be the limits of $\{v_q\}, \{w_q\}$ respectively. Furthermore, interpolation between \eqref{e:vwq+1--1} and \eqref{e:vwq+1--2} gives
\begin{align*}
\|v_{q+1}-v_q\|_{1+\beta}&\leq\|v_{q+1}-v_q\|_{1}^{1-\beta}\|v_{q+1}-v_q\|_{2}^\beta\\
&\leq K_3\delta_{q+1}^{1/2}\lambda_{q+1}^\beta=K_3a^{b^{q+1}(-1+2cb\beta)/2},\\
\|w_{q+1}-w_q\|_{1+\beta}&\leq\|w_{q+1}-w_q\|_{1}^{1-\beta}\|w_{q+1}-w_q\|_{2}^\beta\\
&\leq K_3\delta_{q+1}^{1/2}\lambda_{q+1}^\beta=K_3a^{b^{q+1}(-1+2cb\beta)/2},
\end{align*}
with some constant $K_3$ depending on $K$. Using \eqref{e:delta-q-lambda-q} we can see that if $\beta< \frac{1}{2bc}$ then the exponent in the above estimates are negative. Therefore,  $v\in C^{1, \beta}(\bar{\Omega}), w\in C^{1, \beta}(\bar{\Omega}).$
Moreover, from
\begin{align*}
\|A-\delta_{q+1}\textrm{Id}-\frac{1}{2}\nabla v_q\otimes\nabla v_q-\textrm{sym }\nabla w_q\|_\alpha\leq\sigma_0\delta_{q+1}\rightarrow0, \text{ as } q\rightarrow\infty,
\end{align*}
we have
\begin{align*}
&A=\frac{1}{2}\nabla v\otimes\nabla v-\textrm{sym }\nabla w,\\
&\Det \nabla^2 v=-\textrm{curl }\textrm{curl }A=f.
\end{align*}
By \eqref{e:vwq+1--0}, we get
\begin{align*}
\|v-v_0\|_0\leq\sum_{q\geq1}\delta_{q}^{1/2}\lambda_q^{-\gamma}\leq a^{-\gamma bc-b/2}\leq\frac{\eps}{2},
\end{align*}
provided taking $a$ larger, thus with \eqref{e:first-difference} we have
\begin{align*}
\|v-v^\flat\|_0\leq\eps,
\end{align*}

Now it remains to show that $\beta$ can be taken close to $\frac{1}{5}$ and $p$ can be taken close to $\frac{5}{4},$ which can be obtained from analysing conditions on $b, c, \kappa$ from  Lemma \ref{l:b-c-kappa}. In fact, combining with \eqref{e:l}, \eqref{e:bc-admissible} is equivalent to
\begin{align*}
\delta_{q+2}^2\lambda_{q+1}^{1-2\alpha}>
\delta_{q+1}^{2-\frac{1}{2-\alpha}}\delta_{q}^{\frac{1}{2-\alpha}}\lambda_q^{\frac{2}{2-\alpha}},
\end{align*}
provided $a$ large enough. Inserting \eqref{e:delta-q-lambda-q} and taking logarithm in base $a$ implies
\begin{align*}
(c(1-2\alpha)-2)b^2>\left(\frac{1+2c}{2-\alpha}-2\right)b-\frac{1}{2-\alpha}.
\end{align*}
Thus we require
\begin{align*}
cb((2-\alpha)(1-2\alpha)b-2)>2(2-\alpha)b^2+(1-2(2-\alpha))b-1,
\end{align*}
and $a\gg1$ larger, which then follows from
\begin{equation}\label{e:bc-condition}
b>\frac{2}{(2-\alpha)(1-2\alpha)},\quad c>\frac{2(2-\alpha)b^2-(3-2\alpha)b-1}{b((2-\alpha)(1-2\alpha)b-2)},
\end{equation}
Hence, considering  \eqref{e:bccondition2} and \eqref{e:bc-condition}, taking $\alpha$ arbitrarily small, we have
\begin{align*}
b>1, \quad c>\frac{4b^2-3b-1}{2b(b-1)}=2+\frac{1}{2b}, \quad c>2+\frac{1}{2b}.
\end{align*}
Thus, with $\alpha\to 0$ we can choose $b\to 1$ and $c\to 5/2$, which leads to
\begin{align*}
\beta<\frac{1}{2bc}\to\frac{1}{5}.
\end{align*}
On the other hand, using \eqref{e:delta-q-lambda-q}, it is easy to see that \eqref{e:kappa-admissible} is equivalent to
\begin{align*}
-\kappa b+c\alpha(2-\alpha)b^2<-(2-\alpha)b^2-\kappa+2cb\kappa,
\end{align*}
which can be obtained by taking $a$ sufficient large and $\kappa$ as
\begin{equation}\label{e:kappa-condition}
\kappa>\frac{cb^2\alpha(2-\alpha)+b^2(2-\alpha)}{2cb-1+b},
\end{equation}
Therefore, from \eqref{e:kappa-condition},
\begin{align*}
\kappa>\frac{2b^2}{2cb-1+b}\rightarrow \frac{2}{5}, \text{ as } b\to 1,\,c\to 5/2.
\end{align*}
Thus $$p=\frac{2}{2-\kappa}\rightarrow\frac{5}{4} \text{ as } b\to 1,\, c\to 5/2.$$
We finally gain the proof.

\section{Proof of Iteration Proposition \ref{p:stage}}\label{s:iterate-prop}

In this section we will verify the iteration Proposition \ref{p:stage}, which is a modification of Proposition 5.2 in \cite{Lewicka:2015gz} and  parallel to Proposition 1.1 in \cite{DeLellis:2015wm}.  The difference here is that  we apply Proposition \ref{p:diag} to diagonalise the deficit matrix.  We divide the proof into three steps. In the following  proof, $C$ denotes constants independent of $\alpha, a, b, c, K $ may depend on $A$ and varies from line to line.

\textbf{Step 1. Regularization and diagonalisation.} We first regularize $v_q, w_q, A, D_q$ on length scale $\ell$ to get  $\tilde{v}, \tilde{w}, \tilde{A} $ and $\tilde{D_q}$ respectively by mollifier $\varphi_\ell$. Immediately, from \eqref{e:dq}, \eqref{e:vwq} and Lemma \ref{l:mollification} we obtain
\begin{equation}\label{e:vwtilde}
\begin{split}
\|\tilde{v}-v_q\|_j+\|\tilde{w}-w_q\|_j&\leq C(\|v_q\|_2+\|w_q\|_2)\ell^{2-j}\\
&\leq C(K)\delta_q^{1/2}\lambda_q\ell^{2-j}, j=0, 1,\\
\|\tilde{v}\|_{2+j}+\|\tilde w\|_{2+j}&\leq C(K,j)\delta_q^{1/2}\lambda_q\ell^{-j}, j\geq 0, \\
\end{split}
\end{equation}
and
\begin{equation}
\label{e:aq}
\|\tilde{A}-A\|_j\leq C(j)\|A\|_\kappa\ell^{\kappa-j}, j\geq 0.
\end{equation}
Denote
$$\mathfrak{D}_q=\tilde{A}-\frac{1}{2}\nabla\tilde{v}\otimes\nabla\tilde{v}-\textrm{sym}\nabla\tilde{w}-\delta_{q+2}\textrm{Id},$$
then again by Lemma \ref{l:mollification} and using \eqref{e:vwtilde},  we get
\begin{align*}
\|\mathfrak{D}_q-\delta_{q+1}\textrm{Id}\|_\alpha=&\|\tilde{D_q}-\delta_{q+2}\textrm{Id}+\frac{1}{2}\left((\nabla v_q\otimes\nabla v_q)*\varphi_\ell-\nabla\tilde{v}\otimes\nabla\tilde{v}\right)\|_\alpha\\
\leq&\sigma_0\delta_{q+1}+\delta_{q+2}+CK^2\delta_q\lambda_q^2\ell^{2-\alpha}\\
\leq&3\sigma_0\delta_{q+1},
\end{align*}
provided $K_1>K$ larger, where we have also used
\begin{equation*}
\|\tilde{D_q}\|_\alpha\leq \|D_q\|_\alpha\leq \sigma_0\delta_{q+1}.
\end{equation*}
Hence
\begin{equation}\label{e:dificitq}
\left\|\frac{\mathfrak{D}_q}{\delta_{q+1}}-\textrm{Id}\right\|_\alpha\leq\sigma_1, \quad  \|\mathfrak{D}_q\|_{j+\alpha}\leq C(\sigma_0)\delta_{q+1}\ell^{-j}, j\geq1.
\end{equation}
A direct application of Proposition \ref{p:diag}  yields that there exists $\Phi(x)$ and $\bar{d}(x)>0$ such that
\begin{align*}
\frac{\mathfrak{D}_q}{\delta_{q+1}}+ \textrm{sym }\nabla\Phi=\bar{d}^2 \textrm{Id}.
\end{align*}
Define
$$w_*(x)=-\delta_{q+1}\Phi(x), \quad d(x)=\delta_{q+1}^{1/2}\bar{d}(x),$$
then
$$\mathfrak{D}_q-\textrm{sym }\nabla w_*=d^2\textrm{Id}.$$
From \eqref{e:dificitq} and conclusions of Proposition \ref{p:diag}, we also have
\begin{align}
&\|d\|_{j+\alpha}\leq \delta_{q+1}^{1/2}\|\bar{d}\|_{j+\alpha}\leq C(\sigma_0,j)\delta_{q+1}^{1/2}\ell^{-j}, j\geq 0; \label{e:d-esimate} \\
&\|w_*\|_{j}\leq \delta_{q+1}\|\Phi\|_{j}\leq C(\sigma_0)\delta_{q+1},  j=0, 1;\label{e:w*estimate-1}\\
&\|w_*\|_{1+j}\leq C(\sigma_0,j)\delta_{q+1}\ell^{-j}, j\geq 0.  \label{e:w*estimate-2}
\end{align}

\textbf{Step 2. Adding the first deficit tensor.}
Similar to \cite{Lewicka:2015gz}, applying Lemma \ref{l:Gamma}, we add the first tensor $d^2 e_1\otimes e_1$ as follows.
\begin{align*}
&\bar{v}(x)=\tilde{v}(x)+\frac{1}{\mu}\Gamma_1(d(x), \mu x\cdot e_1),\\
&\bar{w}(x)=\tilde{w}(x)+w_*(x)-\frac{1}{\mu}\Gamma_1(d(x), \mu x\cdot e_1)\nabla\tilde{v}(x)+\frac{1}{\mu}\Gamma_2(d(x), \mu x\cdot e_1)e_1
\end{align*}
with  $\Gamma_i(s, t), i=1, 2$ in \eqref{e:gamma}.
Observe that
\begin{align}
&\|\bar{v}-\tilde{v}\|_j\leq \frac{1}{\mu}\|\Gamma_1\|_j,\label{e:vbarj}\\
&\|\bar{w}-\tilde{w}-w_*\|_j\leq C(j)\frac{1}{\mu}(\|\Gamma_1\|_j\|\tilde{v}\|_1+\|\Gamma_1\|_0\|\tilde{v}\|_{j+1}+\|\Gamma_2\|_j),\label{e:wbarj}
\end{align}
for $ j\geq 0$  So we need to estimate $\|\Gamma_i\|_j $ for $i=1,2$ and $j=0, 1, 2.$  Here $\|\Gamma_i\|_j$ denotes $C^j$ norms of function $x\rightarrow\Gamma_i(d(x), \mu x\cdot e_1)$ and same as  $\|\partial_t\Gamma_i\|_j, \|\partial_s\Gamma_i\|_j.$ Indeed,  using Lemma \ref{l:Gamma} we deduce
\begin{equation}\label{e:gamma1-cjnorm}
\begin{split}
&\|\Gamma_1\|_0+\|\partial_t\Gamma_1\|_0+\|\partial_t^2\Gamma_1\|_0\leq C\|d\|_0\leq C(\sigma_0)\delta_{q+1}^{1/2},\\
&\|\Gamma_1\|_1+\|\partial_t\Gamma_1\|_1\leq C(\|d\|_1+\mu\|d\|_0)\leq C(\sigma_0)\delta_{q+1}^{1/2}\mu,\\
&\|\partial_s\Gamma_1\|_0+\|\partial_s^2\Gamma_1\|_0\leq C,\\
& \|\partial_s\Gamma_1\|_1\leq C\mu,\\
&\|\Gamma_1\|_{1+j}\leq C(\sigma_0,j)\delta_{q+1}^{1/2}\mu^{1+j},\quad j\geq 0,
\end{split}
\end{equation}
where we have used $\mu>\ell^{-1}$ in Lemma \ref{l:b-c-kappa} and \eqref{e:d-esimate} and for the last estimate used the chain rule as in \cite{Conti:2012wl} (alternatively, the specific choice in \eqref{e:gamma} requires only the product rule). Similarly, we have
\begin{equation}\label{e:gamma2-cjnorm}
\begin{split}
&\|\Gamma_2\|_0+\|\partial_t\Gamma_2\|_0+\|\partial_t^2\Gamma_2\|_0\leq C\|d\|_0^2\leq C(\sigma_0)\delta_{q+1},\\
&\|\Gamma_2\|_1+\|\partial_t\Gamma_2\|_1\leq C(\|d^2\|_1+\mu\|d^2\|_0)\leq C(\sigma_0)\delta_{q+1}\mu,\\
&\|\partial_s\Gamma_2\|_0\leq C(\sigma_0)\delta_{q+1}^{1/2},  \quad \|\partial_s^2\Gamma_2\|_0\leq C, \\
&\|\partial_s\Gamma_2\|_1\leq C(\|d\|_1+\mu\|d\|_0)\leq C(\sigma_0)\delta_{q+1}^{1/2}\mu,\\
&\|\Gamma_2\|_{1+j}\leq C(\sigma_0)\delta_{q+1}\mu^{1+j}.
\end{split}
\end{equation}
Finally from \eqref{e:vbarj}, \eqref{e:gamma1-cjnorm} and \eqref{e:vwtilde} we obtain
\begin{equation}\label{e:vbarestimate}
\begin{split}
\|\bar{v}-v_q\|_j&\leq\|\tilde{v}-v_q\|_j+\|\bar{v}-\tilde{v}\|_j\leq C(\sigma_0)\delta_{q+1}^{1/2}\mu^{j-1}, j=0, 1, 2,\\
\|\bar{v}\|_{2+j}&\leq \|\tilde{v}\|_{2+j}+\|\bar{v}-\tilde{v}\|_{2+j}\leq C(\sigma_0,j)\delta_{q+1}^{1/2}\mu^{j+1}, j\geq 0,
\end{split}
\end{equation}
where we have used \eqref{e:bc-admissible}. By \eqref{e:vwtilde}, \eqref{e:wbarj}, \eqref{e:gamma1-cjnorm}, \eqref{e:gamma2-cjnorm}, we also have
\begin{align*}
\|\bar{w}-\tilde{w}-w_*\|_0
\leq\frac{C(\sigma_0)}{\mu}(\delta_{q+1}^{1/2}+\delta_{q+1})\leq C(\sigma_0)\delta_{q+1}^{1/2}\mu^{-1},
\end{align*}
where we have used $\|\nabla\tilde{v}\|_0\leq\|\nabla v\|_0\leq \|A\|_0.$ Moreover,
\begin{equation*}
\begin{split}
\|\bar{w}-\tilde{w}-w_*\|_1&\leq\frac{C(\sigma_0)}{\mu}(\delta_{q+1}^{1/2}\ell^{-1}+\delta_{q+1}^{1/2}\mu+K\delta_{q}^{1/2}\lambda_q\delta_{q+1}^{1/2}+\delta_{q+1}\ell^{-1}+\delta_{q+1}\mu)\\
&\leq C(\sigma_0)\delta_{q+1}^{1/2},\\
\|\bar{w}-\tilde{w}-w_*\|_2
&\leq\frac{C(\sigma_0)}{\mu}(\delta_{q+1}^{1/2}\ell^{-2}+\delta_{q+1}^{1/2}\mu^2+K\delta_{q}^{1/2}\lambda_q\delta_{q+1}^{1/2}\ell^{-1}+\delta_{q+1}\ell^{-2}+\delta_{q+1}\mu^2)\\
&\leq C(\sigma_0)\delta_{q+1}^{1/2}\mu,
\end{split}
\end{equation*}
by appropriately taking $K_2$ large.
Consequently, using \eqref{e:w*estimate-1} and \eqref{e:w*estimate-2} we have
\begin{equation}\label{e:wbarestimate}
\begin{split}
&\|\bar{w}-w_q\|_0\leq\|\tilde{w}-w_q\|_0+\|w_*\|_0+\|\bar{w}-\tilde{w}-w_*\|_0\leq C(\sigma_0)\delta_{q+1},\\
&\|\bar{w}-w_q\|_1\leq\|\tilde{w}-w_q\|_1+\|w_*\|_1+\|\bar{w}-\tilde{w}-w_*\|_1\leq C(\sigma_0)\delta_{q+1}^{1/2},\\
&\|\bar{w}-w_q\|_2\leq\|\tilde{w}-w_q\|_2+\|w_*\|_2+\|\bar{w}-\tilde{w}-w_*\|_2\leq C(\sigma_0)\delta_{q+1}^{1/2}\mu,
\end{split}
\end{equation}
where we have used $\delta_{q+1}^{-1/2}\leq\lambda_q\leq\mu.$
Next we shall estimate $C^j, j=0, 1$ norms of the error matrix $\mathcal{E}_1$ with
\begin{equation*}
\mathcal{E}_1:=\left(\frac{1}{2}\nabla\bar{v}\otimes\nabla\bar{v}+\textrm{sym}\nabla\bar{w}\right)-\left(d^2 e_1\otimes e_1+\frac{1}{2}\nabla\tilde{v}\otimes\nabla\tilde{v}+\textrm{sym}\nabla(\tilde{w}+w_*)\right).
\end{equation*}
In fact, taking gradient of $\bar{v}, \bar{w}$ implies
\begin{align*}
\nabla \bar v&=\nabla \tilde v+\frac{1}{\mu}\partial_s\Gamma_1\nabla d+\partial_t\Gamma_1 e_1,\\
\nabla \bar w&=\nabla(\tilde w+w_*)-\frac{1}{\mu}\partial_s\Gamma_1\nabla\tilde{v}\otimes\nabla d-\partial_t\Gamma_1\nabla\tilde{v}\otimes e_1-\frac{1}{\mu}\Gamma_1\nabla^2\tilde{v}\\
&\quad +\frac{1}{\mu}\partial_s\Gamma_2e_1\otimes\nabla d+\partial_t\Gamma_2e_1\otimes e_1.
\end{align*}
 Direct calculation and using \eqref{e:gamma12} contribute to
 \begin{align*}
 \mathcal{E}_1=&\frac{1}{\mu}\big[\partial_t\Gamma_1\partial_s\Gamma_1\textrm{sym}(e_1\otimes\nabla d)-\Gamma_1\nabla^2\tilde{v}+\partial_s\Gamma_2\textrm{sym}(e_1\otimes\nabla d)\big]\\
 &+\frac{1}{2\mu^2}(\partial_s\Gamma_1)^2\nabla d\otimes\nabla d.
 \end{align*}
\eqref{e:gamma1-cjnorm}, \eqref{e:gamma2-cjnorm} and interpolation inequalities of $C^j$ norms are then  utilized to get
\begin{align*}
 \|\mathcal{E}_1\|_0\leq&\frac{1}{\mu}(\|\partial_t\Gamma_1\|_0\|\partial_s\Gamma_1\|_0\|\nabla d\|_0+\|\Gamma_1\|_0\|\nabla^2\tilde{v}\|_0+\|\partial_s\Gamma_2\|_0\|\nabla d\|_0)\\
 &+\frac{1}{2\mu^2}\|\partial_s\Gamma_1\|_0^2\|\nabla d\|_0^2\\
 \leq&\frac{C(K, \sigma_0)}{\mu}(\delta_{q+1}^{1/2}\delta_{q+1}^{1/2}\ell^{-1}+\delta_{q+1}^{1/2}\delta_q^{1/2}\lambda_q+\delta_{q+1}\ell^{-1})+\frac{C(K, \sigma_0)}{\mu^2}\delta_{q+1}\ell^{-2}\\
 \leq &C(K, \sigma_0)\delta_{q+1}(\mu\ell)^{-1},
 \end{align*}
and
\begin{align*}
 \|\nabla\mathcal{E}_1\|_0\leq&\frac{1}{\mu}(\|\partial_t\Gamma_1\|_1\|\partial_s\Gamma_1\|_0\|\nabla d\|_0+\|\partial_t\Gamma_1\|_0\|\partial_s\Gamma_1\|_1\|\nabla d\|_0\\
 &+\|\partial_t\Gamma_1\|_0\|\partial_s\Gamma_1\|_0\|\nabla d\|_1
 +\|\Gamma_1\|_1\|\nabla^2\tilde{v}\|_0+\|\Gamma_1\|_0\|\nabla^2\tilde{v}\|_1\\
 &+\|\partial_s\Gamma_2\|_1\|\nabla d\|_0+ \|\partial_s\Gamma_2\|_0\|\nabla d\|_1)\\
 &+\frac{1}{\mu^2}(\|\partial_s\Gamma_1\|_0\|\partial_s\Gamma_1\|_1\|\nabla d\|_0^2+\|(\partial_s\Gamma_1)^2\|_0\|\nabla d\|_0\|\nabla d\|_1)\\
 \leq&\frac{C(K, \sigma_0)}{\mu}(\delta_{q+1}^{1/2}\mu \delta_{q+1}^{1/2}\ell^{-1}+\delta_{q+1}^{1/2}\mu\delta_q^{1/2}\lambda_q+\delta_{q+1}^{1/2}\delta_q^{1/2}\lambda_q\ell^{-1}+\delta_{q+1}\ell^{-2})\\
&+\frac{C(K, \sigma_0)}{\mu^2}(\delta_{q+1}\mu\ell^{-2}+\delta_{q+1}\ell^{-3})\\
\leq& C(K, \sigma_0)\delta_{q+1}\ell^{-1}.
 \end{align*}
Thus we get
\begin{equation}
\label{e:error-1}
\|\mathcal{E}_1\|_0\leq C(K, \sigma_0)\delta_{q+1}(\mu\ell)^{-1}, \quad \|\nabla\mathcal{E}_1\|_0\leq C(K, \sigma_0)\delta_{q+1}\ell^{-1}.
\end{equation}

 \textbf{Step 3. Adding the second deficit tensor and conclusion.} Similar to Step 2, to add tensor $d^2e_2\otimes e_2,$  we construct our final $v_{q+1}, w_{q+1}$ through the following
\begin{align*}
&v_{q+1}(x)=\bar{v}(x)+\frac{1}{\lambda_{q+1}}\Gamma_1(d(x), \lambda_{q+1} x\cdot e_2),\\
&w_{q+1}(x)=\bar{w}(x)-\frac{1}{\lambda_{q+1}}\Gamma_1(d(x), \lambda_{q+1} x\cdot e_2)\nabla \bar{v}(x)+\frac{1}{\lambda_{q+1}}\Gamma_2(d(x), \lambda_{q+1} x\cdot e_2)e_2.
\end{align*}
Parallel to Step 2, a similar calculation will contribute to
\begin{align*}
&\|\Gamma_1\|_0+\|\partial_t\Gamma_1\|_0+\|\partial_t^2\Gamma_1\|_0\leq C(\sigma_0)\delta_{q+1}^{1/2},\\
&\|\Gamma_1\|_1+\|\partial_t\Gamma_1\|_1\leq C(\sigma_0)\delta_{q+1}^{1/2}\lambda_{q+1},\\
&\|\partial_s\Gamma_1\|_0+\|\partial_s^2\Gamma_1\|_0\leq C,\\
&\|\partial_s\Gamma_1\|_1\leq C\lambda_{q+1}, \quad \|\Gamma_1\|_2\leq C(\sigma_0)\delta_{q+1}^{1/2}\lambda_{q+1}^2.
\end{align*}
and
\begin{align*}
&\|\Gamma_2\|_0+\|\partial_t\Gamma_2\|_0+\|\partial_t^2\Gamma_2\|_0\leq C(\sigma_0)\delta_{q+1},\\
&\|\Gamma_2\|_1+\|\partial_t\Gamma_2\|_1\leq C(\sigma_0)\delta_{q+1}\lambda_{q+1},\\
&\|\partial_s\Gamma_2\|_0\leq C(\sigma_0)\delta_{q+1}^{1/2},  \quad \|\partial_s^2\Gamma_2\|_0\leq C, \\
&\|\partial_s\Gamma_2\|_1\leq C(\sigma_0)\delta_{q+1}^{1/2}\lambda_{q+1},\quad
\|\Gamma_2\|_2\leq C(\sigma_0)\delta_{q+1}\lambda_{q+1}^2.
\end{align*}
Here $\|\Gamma_i\|_j$ denotes $C^j$ norms of function $x\rightarrow\Gamma_i(d(x), \lambda_{q+1} x\cdot e_2)$ as before and same as $\|\partial_t\Gamma_i\|_j, \|\partial_s\Gamma_i\|_j.$
Thus with the formulae of $v_{q+1}$ and $w_{q+1}$, direct calculation gives us for $j=0, 1, 2,$
\begin{equation}\label{e:vwq+1}
\|\bar{v}-v_{q+1}\|_j+\|\bar{w}-w_{q+1}\|_j\leq C(\sigma_0)\delta_{q+1}^{1/2}\lambda_{q+1}^{j-1}.
\end{equation}
Summing up \eqref{e:vbarestimate}, \eqref{e:wbarestimate} and \eqref{e:vwq+1}, we have
\begin{align*}
&\|v_{q+1}-v_q\|_0\leq \|v_{q+1}-\bar{v}\|_0+\|\bar{v}-v_q\|_0\leq C(\sigma_0)\delta_{q+1}^{1/2}\mu^{-1}\\
&\|v_{q+1}-v_q\|_1\leq  \|v_{q+1}-\bar{v}\|_1+\|\bar{v}-v_q\|_1\leq C(\sigma_0)\delta_{q+1}^{1/2},\\
&\|v_{q+1}-v_q\|_2\leq C(\sigma_0)\delta_{q+1}^{1/2}(\mu+\lambda_{q+1})\leq K\delta_{q+1}^{1/2}\lambda_{q+1},
\end{align*}
and
\begin{align*}
&\|w_{q+1}-w_q\|_0\leq \|w_{q+1}-\bar{w}\|_0+\|\bar{w}-w_q\|_0\leq C(\sigma_0)\delta_{q+1},\\
&\|w_{q+1}-w_q\|_1\leq \|w_{q+1}-\bar{w}\|_1+\|\bar{w}-w_q\|_1\leq C(\sigma_0)\delta_{q+1}^{1/2},\\
&\|w_{q+1}-w_q\|_2\leq C(\sigma_0)\delta_{q+1}^{1/2}(\mu+\lambda_{q+1})\leq K\delta_{q+1}^{1/2}\lambda_{q+1},
\end{align*}
which imply \eqref{e:vwq+1--0} and \eqref{e:vwq+1--1} for some constant $M$ depending on $\sigma_0.$ Since $\delta_{q}^{1/2}\lambda_q\leq\delta_{q+1}^{1/2}\lambda_{q+1}$ by Lemma \ref{l:b-c-kappa}, we also arrive at \eqref{e:vwq+1--2}. Moreover, for the second matrix error
 \begin{align*}
 \mathcal{E}_2&:=\left(\frac{1}{2}\nabla v_{q+1}\otimes\nabla v_{q+1}+\textrm{sym}\nabla w_{q+1}\right)-\left(d^2 e_2\otimes e_2+\frac{1}{2}\nabla\bar{v}\otimes\nabla\bar{v}+\textrm{sym}\nabla\bar{w}\right)\\
 &=\frac{1}{\lambda_{q+1}}\big[\partial_t\Gamma_1\partial_s\Gamma_1\textrm{sym}(e_2\otimes\nabla d)-\Gamma_1\nabla^2\bar{v}+\partial_s\Gamma_2\textrm{sym}(e_2\otimes\nabla d)\big]\\
 &\quad +\frac{1}{2\lambda_{q+1}^2}(\partial_s\Gamma_1)^2\nabla d\otimes\nabla d,
 \end{align*}
in a similar way to Step 2 we are able to deduce
 \begin{equation}\label{e:error-2}
 \|\mathcal{E}_2\|_0\leq C(K, \sigma_0)\delta_{q+1}\mu\lambda_{q+1}^{-1}, \qquad
 \|\nabla\mathcal{E}_2\|_0\leq C(K, \sigma_0)\delta_{q+1}\mu.
 \end{equation}
Since
\begin{align*}
D_{q+1}=&A_{q+1}-\frac{1}{2}\nabla v_{q+1}\otimes\nabla v_{q+1}-\textrm{sym}\nabla w_{q+1}\\
=&A_{q+1}-(\frac{1}{2}\nabla\tilde{v}\otimes\nabla\tilde{v}+\textrm{sym}\nabla\tilde{w}+\textrm{sym}\nabla w_*+d^2\textrm{Id}+\mathcal{E}_1+\mathcal{E}_2)\\
=&A-(\delta_{q+2}\textrm{Id}+\frac{1}{2}\nabla\tilde{v}\otimes\nabla\tilde{v}+\textrm{sym}\nabla\tilde{w}+\mathfrak{D}_q+\mathcal{E}_1+\mathcal{E}_2)\\
=&A-\tilde{A}-\mathcal{E}_1-\mathcal{E}_2,
\end{align*}
by \eqref{e:aq}, \eqref{e:error-1} and \eqref{e:error-2}, we have
\begin{align*}
\|D_{q+1}\|_0\leq&\|A-\tilde{A}\|_0+\|\mathcal{E}_1\|_0+\|\mathcal{E}_2\|_0\\
\leq &C\ell^\kappa+C(K, \sigma_0)\delta_{q+1}(\mu\ell)^{-1}+C(K, \sigma_0)\delta_{q+1}\mu\lambda_{q+1}^{-1}\\
\leq&\frac{\sigma_0}{3}\delta_{q+2}\lambda_{q+1}^{-\alpha},
\end{align*}
by Lemma \ref{l:b-c-kappa} and taking $a$ large enough.
As for $\|\nabla D_{q+1}\|_0$, it easily follows
\begin{align*}
\|\nabla D_{q+1}\|_0\leq&\|A-\tilde{A}\|_1+\|\nabla\mathcal{E}_1\|_0+\|\nabla\mathcal{E}_2\|_0\\
\leq &C\ell^{\kappa-1}+C(K, \sigma_0)\delta_{q+1}\ell^{-1}+C(K, \sigma_0)\delta_{q+1}\mu\\
\leq&\frac{\sigma_0}{3}\delta_{q+2}\lambda_{q+1}^{1-\alpha},
\end{align*}
provided $a$ large enough. Thus we gain \eqref{e:dq+1} and then complete the proof.

\Acknowledgements{The authors would like to thank the hospitality of the Max-Plank Institute of
Mathematics in the Sciences, and gratefully acknowledge the support of the ERC
Grant Agreement No. 724298.}


\end{document}